\documentclass{amsart}

\usepackage{amsthm,amsfonts,amsmath,amssymb,latexsym,epsfig,eucal,ae,upref}
\usepackage{undertilde}

\newtheorem{theorem}{Theorem}
\newtheorem{lemma}[theorem]{Lemma}

\newtheorem{corollary}[theorem]{Corollary}

\newenvironment{remark}{\medskip \refstepcounter{theorem}
\noindent  {\bf Remark \thetheorem}.\rm}{\,}

\def\BOne{{\mathchoice {\rm 1\mskip-4mu l} {\rm 1\mskip-4mu l}
		  {\rm 1\mskip-4.5mu l} {\rm 1\mskip-5mu l}}}

\def\tg{{\tilde g}}

\def\<{\langle}

\def\n{\nabla}
\def\>{\rangle}

\def\a{\alpha}

\def\tm{\tilde{M}}
\def\tn{\tilde{n}}
\def\tg{\tilde{g}}
\def\lat{\utilde{A}} 
\def\mb#1{{\mathbb #1}}

\begin{document}

\title[Low energy canonical immersions]{Low energy canonical immersions into 
hyperbolic manifolds and standard spheres}
\author[H. del Rio]{Heberto del Rio}
\address{Department of Mathematics \& Computer Sciences, 
Barry University, 11300 Northeast Second Avenue,
Miami Shores, FL 33161, U.S.A.}
\author[W. Santos]{Walcy Santos}
\address{Instituto de Matem\'atica, Universidade Federal do Rio de Janeiro,
68530, 21941-909, Rio de Janeiro, BRAZIL.}
\author[S.R. Simanca]{Santiago R. Simanca}
\address{Department of Mathematics, University of Miami, Coral Gables, FL 33124, U.S.A.}
\email{hdelrio@mail.barry.edu, walcy@im.ufrj.br, srsimanca@gmail.com}

\begin{abstract}
We consider critical points of the functionals $\Pi$ and $\Psi$ defined as
the global $L^2$-norm of the second fundamental form and mean curvature vector
of isometric immersions of compact Riemannian manifolds into a background Riemannian
manifold, respectively, as functionals over the space of deformations of the 
immersion. We prove gap theorems for these functionals into hyperbolic manifolds,
and show that the celebrated gap theorem for minimal immersions into the 
sphere can be cast as a theorem about critical points of these functionals of constant
mean curvature function, and whose second fundamental form is suitably small in 
relation to it. In this case, the various type of minimal submanifolds that can occur 
at the pointwise upper bound on the norm of the second fundamental form are realized 
by manifolds of nonnegative Ricci curvature, and of these, the Einstein ones are 
distinguished from the others by being those that are immersed on the sphere as 
critical points of $\Pi$. 
\end{abstract}

\maketitle
\section{Introduction}
Let $n$ and $\tn$ be any two integers such that $2\leq n \leq \tn-1$. 
Let $M$ be a closed Riemannian manifold of dimension $n$ 
isometrically immersed into a background Riemannian manifold $(\tm,\tg)$
of dimension $\tn$. Thus, we have an immersion $f:M \rightarrow \tm$ of $M$ 
into $\tm$ and the Riemannian metric on $M$ is that induced by the ambient 
space metric $\tg$. We let $\a$ and $H$ be the second fundamental form and 
mean curvature vector of the immersion, respectively. If $d\mu=d\mu_M$ denotes
the Riemannian measure on $M$, we define 
\begin{equation}
\Pi(M)= \int_{M} \| \a \|^2 \, d\mu \, , \label{sf}
\end{equation}
\begin{equation}
\Psi(M)= \int_{M} \| H \|^2 \, d\mu \, , \label{mc}
\end{equation}
and view them as functionals defined over the space of all isometric 
immersions of $M$ into $\tm$. $\Psi$ is generally named the Willmore 
functional in the 
literature. Willmore used it to study surfaces in $\mb{R}^3$ \cite{wi}. 
This functional had
been studied earlier by Blaschke \cite{bla} and Thomsen \cite{th} also.
 
We shall say that a Riemannian manifold $(M,g)$ is {\it canonically placed} 
into $(\tm,\tg)$ if it admits an isometric 
immersion into the latter that is a critical point of $\Pi$. Thus, we have arbitrary families 
of metrics $g_t$ in $M$ that deform $g$, $g_0=g$, and that admit isometric immersions 
$f_t: M \hookrightarrow \tm$.
Then, $(M,g)$ is canonically placed into $(\tm,\tg)$ if $f_0$ is a critical point of
$\Pi_t$. The use of the functional $\Psi_t$ instead leads to a related, and 
alternative, notion of a canonical placing of $(M,g)$ into $(\tm,\tg)$. 

In the case when the background manifold is the space form of curvature 
$c$, we 
introduce for consideration a third functional given by 
\begin{equation}\label{nf}
\Theta_c (M)=\int_{M} (n(n-1)c+\| H\|^2)\,  d\mu \, . 
\end{equation}
When $c=1$, we denote $\Theta_c$ simply by $\Theta$.

The main purpose of our work is the study of the critical points of lowest energy of $\Pi$ 
and $\Psi$ when the background manifold is the sphere $\mb{S}^{n+p}$ with its standard 
metric. Our results reinterpret the celebrated gap theorem of Simons \cite{si} in these
contexts, and show that that result is in effect one about critical points of $\Psi$ and
$\Pi$ of constant mean curvature function, and whose density is pointwise small.

Let us begin our explanation by recalling the gap theorem in a form suitable 
to our work:

\begin{theorem}\label{t1} \cite[Theorem 5.3.2, Corollary 5.3.2]{si} 
\cite[Main Theorem]{cdck} \cite[Corollary 2]{bl} 
Suppose that $M^n\hookrightarrow \mb{S}^{n+p}$ is an isometric minimal
immersion. Assume that the pointwise inequality 
$\| \a \|^2 \leq np/(2p-1)$ holds everywhere. Then 
\begin{enumerate}
\item Either $\| \a\|^2=0$, or 
\item $\| \a \|^2 =np/(2p-1)$ if, and only if, either $p=1$ and $M^n$ is the minimal 
Clifford torus $\mb{S}^{m}(\sqrt{m/n})\times \mb{S}^{n-m}(\sqrt{(n-m)/n})
\subset \mb{S}^{n+1}$, $1\leq m< n$, with $\| \a\|^2=n$, 
or $n=p=2$ and $M$ is the real projective
plane embedded into $\mb{S}^4$ by the Veronese map with $\| \a\|^2=4/3$.
\end{enumerate}
\end{theorem}

Thus, among minimal $n$-manifolds $M\hookrightarrow \mb{S}^{n+p}$, or critical points of 
the volume functional of $M$, the lowest value that 
$\| \a \|^2$ can achieve is isolated and it is achieved by submanifolds inside equatorial
totally geodesic spheres, while its first nonzero value is the constant $pn/(2p-1)$, 
which is achieved by all the Clifford toruses when $p=1$, or by a minimal real 
projective plane of scalar curvature $2/3$ in $\mb{S}^4$. In this note we prove
that this gap theorem is really characterizing compact Riemannian 
manifolds $(M,g)$ 
that are isometrically embedded into $\mb{S}^{n+p}$ as critical points of the functionals 
(\ref{mc}) or (\ref{nf}), which have constant mean curvature function, and 
whose density is pointwise small. The latter conditions make of the said critical points
absolute minimums of the functional (\ref{mc}). The gap theorem then produces 
two type of such critical points, and those that achieve the upper bound for
$\| \alpha\|^2$ carry metrics of nonnegative Ricci tensor. Among them, the ones that
are Einstein are distinguished by being, in addition, critical points of (\ref{sf}).

For an immersion $M \hookrightarrow \tm$, we let $\nu_H$ denote the normal vector in the 
direction of the mean curvature vector $H$, and denote by 
$A_{\nu_H}$ and $\nabla^{\nu}$ the shape operator in the direction of $\nu_H$ and covariant 
derivative of the normal bundle, respectively. 
We consider immersions that satisfy the estimates  
\begin{equation} \label{es}
\begin{array}{rcl}
{\displaystyle -\lambda \| H\|^2 -n} & \leq & 
{\displaystyle {\rm trace}\, A_{\nu_H}^2-\| H\|^2 -\| \nabla^\nu \nu_H \|^2} \vspace{1mm} \\ 
 & \leq & {\displaystyle \| \a \|^2 -\| H\|^2 -\| \nabla^\nu \nu_H \|^2} 
\leq {\displaystyle \frac{np}{2p-1}}
\end{array}
\end{equation} 
for some constant $\lambda$. Notice that $\| A_{\nu_H}\|^2={\rm trace}\, A_{\nu_H}^2$ is 
bounded above by
$\| \alpha \|^2$, and so the second of the inequalities above is always true. 

Our first result is the following:

\begin{theorem}\label{mt}
Suppose that $(M^n,g)$ is a closed Riemannian manifold 
isometrically immersed into $\mb{S}^{n+p}$ as a critical point of the functional
$\Psi$ above, and having constant mean curvature function $\| H\|$. Assume that 
the immersion is such that  {\rm (\ref{es})} holds for some  
constant $\lambda \in [0,1/2)$. Then $M$ is minimal, and
so it is a critical point of the functional $\Theta$ also, 
$0\leq \| \a \|^2 \leq np/(2p-1)$, and either
\begin{enumerate}
\item $\| \a\|^2=0$, in which case $M$ lies in an equatorial sphere, or 
\item $\| \a\|^2=np/(2p-1)$, in which case either 
$p=1$ and $M^n$ is the Clifford torus 
$\mb{S}^{m}(\sqrt{m/n})\times \mb{S}^{n-m}(\sqrt{(n-m)/n})
\subset \mb{S}^{n+1}$, $1\leq m< n$, with $\| \a\|^2=n$, 
or $n=p=2$ and $M^2$ is the real projective
plane embedded into $\mb{S}^4$ by the Veronese map with $\| \a\|^2=4/3$ and
scalar curvature $2/3$, all cases of metrics with nonnegative Ricci tensor.
\end{enumerate}
\end{theorem}

We now consider immersions that satisfy the estimates  
\begin{equation} \label{es2}
\begin{array}{rcl}
{\displaystyle -\lambda \| H\|^2 -1 } & \leq & 
{\displaystyle {\rm trace}\, \left( \frac{1}{\| H\|} A_{\nu_H} \sum A_{\nu_j}^2\right)  
-\frac{1}{2} \| \alpha\|^2 -\| \nabla^\nu \nu_H\|^2-\| H\|^2} \\ 
 & \leq & {\displaystyle \| \a \|^2 -\| H\|^2-\| \nabla^\nu \nu_H\|^2} 
\leq {\displaystyle \frac{np}{2p-1} }
\end{array}
\end{equation} 
for some constant $\lambda$. Here, $\{ \nu_j\}$ is an orthonormal frame of the normal bundle.
We shall see later on that the second of the inequalities above is always true. 

Our second result distinguishes further the critical points of $\Psi$ 
obtained in Theorem \ref{mt}.

\begin{theorem}\label{mt2}
Suppose that $(M^n,g)$ is a closed Riemannian manifold 
isometrically immersed into $\mb{S}^{n+p}$ as a critical point of the functional
$\Pi$ above, and having constant mean curvature function. Assume that the immersion is such that  
{\rm (\ref{es2})} holds for some constant $\lambda \in [0,1)$. 
Then $M$ is minimal, and so a critical point of $\Psi$ and $\Theta$ also, 
$0\leq \| \a \|^2 \leq np/(2p-1)$, and either  
\begin{enumerate}
\item $\| \a\|^2=0$, in which case $M$ lies inside an equatorial sphere, or
\item $\| \a\|^2=np/(2p-1)$, in which case either $n=p=2$ and $M$ is a 
minimal real projective plane with an Einstein metric 
embedded into $\mb{S}^4$, or $p=1$, $n=2m$ and 
$M$ is the Clifford torus $\mb{S}^{m}(\sqrt{1/2})\times \mb{S}^{m}
(\sqrt{1/2})\subset \mb{S}^{n+1}$ with its Einstein product metric.
\end{enumerate}
\end{theorem}

The following observation follows easily, but it emphasizes the fact
that among these submanifolds, we have some that are critical points
of the total scalar curvature functional among metrics in $M$ 
realized by isometric immersions into the sphere $\mb{S}^{n+p}$.

\begin{corollary}\label{co4}
Let $(M,g)$ be a closed Riemannian manifold that is canonically placed in $\mb{S}^{n+p}$ 
with constant mean curvature function and satisfying {\rm (\ref{es2})} for some $\lambda \in
[0,1)$. Then $M$ is a minimal critical point of the total scalar curvature functional under 
deformations of the isometric immersion, and either $\| \a \|^2=0$ 
or $\| \a \|^2=np/(2p-1)$. In the latter case, $(M,g)$ is Einstein and the two possible 
surface cases in codimension $p=1$ and $p=2$ correspond to Einstein manifolds that are 
associated to different critical values of the total scalar curvature.
\end{corollary}

Thus, among Riemannian manifolds of constant mean curvature function, the pointwise estimates
(\ref{es2}), $\lambda \in [0,1)$, for a critical point of the $L^2$-norm of the second 
fundamental form implies that this is also a minimal critical point of the $L^2$-norm of the 
mean curvature vector, and therefore of the total scalar curvature functional,  
under deformations of the immersion. Of course, this in itself does not imply necessarily that 
$(M,g)$ is Einstein, as the space of metrics on $M$ that can be realized by isometric embeddings 
into the sphere does not have to equal the space of all Riemannian metrics on $M$.  
But if $\| \a\|^2=np/(2p-1)$, the canonically placed submanifold $(M,g)$ is an Einstein 
critical point of the total scalar curvature functional among metrics on $M$ that can be 
realized by isometric immersions into the standard sphere $\mb{S}^{n+p}$. The 
cases in Theorem \ref{t1} that are excluded in Theorem \ref{mt2} correspond to critical points 
of $\Psi$ that are not critical points of $\Pi$. 

That the functional $\Pi$ distinguishes the
symmetric minimal Clifford torus in $\mb{S}^{2m+1}$ from the others have
been observed previously \cite{s7}. Cf. \cite{glw}, \S 3, and \cite{hl}, 
p. 366. It has been observed also that even on
the sphere, the functionals $\Pi$ and $\Psi$ contain critical points that are
not minimal submanifolds \cite{s7}.  

It is slightly easier to prove Theorem \ref{mt} by replacing the role that
the functional $\Psi$ plays by that of the functional $\Theta$, and derive 
the same conclusion. The point is not the use of
critical points of $\Theta$ versus those of $\Psi$. Rather, since the 
curvature of the sphere is positive, if the second fundamental form is 
pointwise small in relation to the mean curvature vector, the constant mean
curvature function condition 
forces the critical points of these functionals to be the same, and minimal. 

The pointwise estimates (\ref{es}) and (\ref{es2}) that we assume on the immersion
quantities are to ensure that the critical submanifolds $(M,g)$ of $\Psi$  
or $\Pi$ that we consider are not too far off from being totally geodesic, 
among immersions with constant mean curvature function. 

We see that the role that the volume functional plays in the gap theorem can
be reinterpreted using the functional $\Psi$, changing accordingly the minimality
condition for having a constant mean curvature function, and the estimates 
(\ref{es}). This singles out compact immersions into the sphere that minimize 
the volume, and so $\Theta$. 
This type of critical points is further distinguished by selecting those 
that are canonically placed in $\mb{S}^{n+p}$ and satisfy the stricter estimates
(\ref{es2}), and among these, selecting some of those
for which the metric $g$ on $M$ is Einstein.  
That the critical points of $\Psi$ and $\Pi$ are this closely interconnected 
is a consequence of Gauss' identity and the fact that the sectional curvature 
of the standard metric on the sphere is the (positive) constant $1$. 

Although we cast Simons' gap theorem in this manner, Theorems \ref{mt} and \ref{mt2} are 
both consequences of nonlinear versions of the first eigenvalue of the Laplacian of the 
metric $g$ on $M$, when $(M,g)$ is a critical point of low density of the functionals in 
question. The natural gap theorem for the functionals $\Psi$ or $\Pi$ themselves that we 
can derive is somewhat dual to that proven for minimal immersions into spheres, and
occur on quotients of space forms of negative curvature instead. 

We recall that a closed hyperbolic manifold is of 
the form $\mb{H}^m/\Gamma$ for $\Gamma$ a torsion-free discrete group of isometries of 
$\mb{H}^m$.  We have the following

\begin{theorem} \label{th5}
Let $M$ be a critical point of {\rm (\ref{mc})} on a hyperbolic compact 
manifold $\mb{H}^{n+p}/\Gamma$. If the pointwise inequality $0\leq \| \alpha \|^2-
\frac{1}{2}\| H\|^2 -\| \nabla^\nu \nu_H\|^2 \leq n$ holds on $M$, then either $\| H\|^2 =0$ and 
$M$ is minimal, or $\| \alpha\|^2=\frac{1}{2}\| H\|^2+n=\| A_{\nu_H}\|^2$ and $M$ is a 
nonminimal submanifold whose mean curvature vector is a covariantly constant section of its 
normal bundle.
\end{theorem}

\begin{theorem} \label{th6}
Let $M$ be a critical point of {\rm (\ref{sf})} on a hyperbolic compact 
manifold $\mb{H}^{n+p}/\Gamma$. If the pointwise inequality 
$$
\| \alpha\|^2 \left( \left( 3-\frac{n}{2}\right) \| \alpha\|^2 -\| H\|^2\right) 
\leq (n\| \alpha\|^2 + 2 \| H \|^2) 
$$
holds, then either $\| H \|^2 =0$ and $M$ is a minimal submanifold, 
or $n\leq 5$, the equality above holds, $\| \alpha \|^2=\| A_{\nu_H}\|^2$, 
and $M$ is a submanifold whose mean curvature vector is a covariantly constant 
section of its normal bundle.
\end{theorem}

It is worth mentioning that though in general we work assuming that $n\geq 2$ in order to 
see the effect of curvature quantities, Theorems \ref{th5} and \ref{th6} remain true for $n=1$, 
case where being minimal and totally geodesic are equivalent notions, and where the two theorems
yield the same result. 

\section{Critical points of the Lagrangians}
Consider a closed Riemannian manifold $(\tm, \tg)$, and let $M$ be an 
submanifold of $\tm$. With the metric $g$ induced by $\tg$, $M$ becomes 
a Riemannian manifold. We denote by $\n^{\tg}$ and $\n^g$ the Levi-Civita
connections of $\tg$ and $g$, respectively, and by
$\a$ the second fundamental form of the isometric immersion. 
The dimensions of $M$ and $\tm$ are
$n$ and $\tn$, respectively.

We have Gauss' identity
\begin{equation}
\n^{\tg}_X Y = \n^{g}_X Y + \a(X,Y)\, . \label{ga}
\end{equation}
If $N$ is a section of the normal bundle $\nu(M)$, 
the shape operator $A_N$ is defined by
$$
A_N X= - \pi_{TM}(\n^{\tg}_X N )\, ,
$$
where in the right side above,
$N$ stands for an extension of the original section to
a neighborhood of $M$.
If $\n^{\nu }$ is the connection on $\nu(M)$ induced by $\n^{\tg}$, we have
Wiengarten's identity
\begin{equation}
\n^{\tg}_X N = -A_N X + \n^{\nu}_X N \, .\label{we}
\end{equation}
For a detailed development of these and some of the expressions that 
follow, see \cite{md}.

Gauss's identity implies Gauss' equation
\begin{equation}
g(R^g(X,Y)Z,W)= \tg(R^{\tg}(X,Y)Z,W)+\tg(\a(X,W),\a(Y,Z))-\tg(\a(X,Z),\a(Y,W))
\, .\label{cu}
\end{equation}
Here, $R^g$ stands for the Riemann curvature tensor of the corresponding 
metric $g$, and $X$, $Y$, $Z$ and $W$ are vector fields in $\tm$ 
tangent to the submanifold $M$. 

Let $\{e_1, \ldots , e_{\tn}\}$ be an orthonormal frame for $\tg$ in a tubular
neighborhood of $M$ such that $\{e_1, \ldots , e_{n}\}$ constitutes an
orthonormal frame for $g$ on points of $M$. We denote by $H$ the mean curvature
vector, the trace of $\a$. The immersion $M \hookrightarrow \tm$ is said to be
minimal if $H=0$. 
By (\ref{cu}), the Ricci tensors $r_g$ and $r_{\tg}$ are related to
each other by the expression 
\begin{equation}\label{rt}
\begin{array}{rcl}
r_g(X,Y) \! \! \! & = & \! \! \! \sum_{i=1}^n \tg(R^{\tg}(e_i,X)Y,e_i)+
\tg(H,\a(X,Y))-\sum_{i=1}^n \tg(\a(e_i,X),\a(e_i,Y)) \\ & = & 
r_{\tg}(X,Y) - \sum_{i=n+1}^{\tn} \tg(R^{\tg}(e_i,X)Y,e_i) +\tg(H,\a(X,Y))-
\\ & & \mbox{} \hspace{2in} \sum_{i=1}^n \tg(\a(e_i,X),\a(e_i,Y)) \, ,
\end{array}
\end{equation}
and the scalar curvatures $s_g$ and $s_{\tg}$ by the expression
\begin{equation}
\begin{array}{rcl}
s_g & = & s_{\tg} -2\sum_{i=1}^n\sum_{j=n+1}^{\tn}K_{\tg}(e_i,e_j) -
K_{\tg}(e_i, e_j)+ \tg(H,H)-\tg(\a,\a) \vspace{1mm} \\
& = & \sum_{i,j\leq n}K_{\tg}(e_i,e_j)+\tg(H,H)- \tg(\a,\a) \, ,\label{sc}
\end{array}
\end{equation}
where $K_{\tg}(e_i,e_j)$ is the $\tg$-curvature of the section spanned by
the orthonormal vectors $e_i$ and $e_j$, and $\tg(H,H)$ and $\tg(\a,\a)$ are
the squared-norms of the mean curvature vector $H$ and the form $\a$,
respectively.

The critical submanifolds for the functionals $\Pi$ and $\Psi$ in
(\ref{sf}) and (\ref{mc}) under deformations of the immersion $f$ are 
described in full generality in \cite[Theorem 3.10]{s7}. We recall the
equations they satisfy adapted to the case of interest here. 
We denote by $S_c^n$ the $n$th dimensional space form of curvature $c$.
For convenience, we use the standard double index summation convention.

\begin{theorem} \cite[Theorem 3.10]{s7} \label{th2}
Let $M$ be an $n$-manifold of codimension $p$ isometrically immersed
into the space form $S^{n+p}_c$. Let $\{ e_1, \ldots, e_n\}$ be an orthonormal 
frame of tangent vectors to $M$, and $\{ \nu_1,\ldots, \nu_p\}$
be an orthonormal frame of the normal bundle of the immersion
such that $H=h\nu_1$. If $\a(e_i,e_j)=h_{ij}^r \nu_r$, then 
$M$ is a critical point of $\Pi$ if, and only if,
$$
2\Delta h  =  2 c h -2h\| \n_{e_i}^{\nu}
\nu_1 \|^2 -
h\| \a\|^2 +2{\rm trace}\, A_{\nu_1}A_{\nu_k}^2 \, ,
$$
and for all $m$ in the range $2\leq m \leq p$, we have that
$$
0  =  2\< e_i(h)\n_{e_i}^{\nu}\nu_1,\nu_m\> +he_i\< \n_{e_i}^{\nu}\nu_1,
\nu_m\> - h\<  \n_{e_i}^{\nu}
\nu_1, \n_{e_i}^{\nu} \nu_m \>
+2{\rm trace}\, A_{\nu_m}A_{\nu_k}^2\, .
$$
$M$ is a a critical point of the
functional {\rm (\ref{mc})} if, and only if,
$$
2\Delta h  =  2cn h 
-2h\| \n_{e_i}^{\nu}\nu_1 \|^2 -
h^3 +2h\, {\rm trace}\, A_{\nu_1}^2 \, ,
$$
and for all $m$ in the range $2\leq m \leq p$, we have that
$$
0 = 4e_i(h)\< \n_{e_i}^{\nu}\nu_1,
\nu_m\>+2he_i \< \n_{e_i}^{\nu}\nu_1,\nu_m\>
- 2h\< \n_{e_i}^{\nu}
\nu_1, \n_{e_i}^{\nu} \nu_m \>
+2h\, {\rm trace}\, A_{\nu_1}A_{\nu_m}\, .
$$
And $M$ is a critical point of the functional {\rm (\ref{nf})} if, and only,
if, the last $p-1$ equations above hold, and the one before these is 
replaced by 
$$
2\Delta h  =  (3n-n^2)ch 
-2h\| \n_{e_i}^{\nu}\nu_1 \|^2 -
h^3 +2h\, {\rm trace}\, A_{\nu_1}^2 \, .
$$
\qed
\end{theorem}

For hypersurfaces into an Einstein background $(\tm,\tg)$, these critical
point equations can be
described completely using only the principal curvatures
$k_1, \ldots, k_n$. Observe that we have 
$h=k_1 + \cdots + k_n$ and $\| \a\|^2 = k_1^2 + 
\cdots + k_n^2$, respectively. 

\begin{theorem} \cite[Theorem 3.3, Corollary 3.4]{s7} \label{th1}
Let $M$ be a hypersurface in an Einstein manifold $(\tm,\tg)$. Assume that
$k_1, \ldots , k_n$ are the principal curvatures, with associated orthonormal 
frame of principal directions $e_1, \ldots, e_n $. Let $\nu$ be a normal
field along $M$. Then $M$ is a critical point of the 
functional {\rm (\ref{sf})} if, and only if, 
$$
2\Delta h = 2 (k_1 K_{\tg}(e_1,\nu)+\cdots + k_n K_{\tg}(e_n,\nu))
 -h\| \a\|^2 +
2(k_1^3 + \cdots + k_n^3) \, ,
$$
and a critical point of the functional {\rm (\ref{mc})} if, and only if, 
$$
2\Delta h = 2h(K_{\tg}(e_1,\nu)+\cdots + K_{\tg}(e_n,\nu))
 +2h\| \a\|^2 -h^3 \, .
$$
In particular, a hypersurface $M$ in
$S_c^{n+1}$ is a critical point for the functional {\rm (\ref{sf})} 
if, and only if, its mean curvature function $h$ satisfies the equation
$$
2\Delta h = 2 ch -h\| \a\|^2 +
2(k_1^3 + \cdots + k_n^3) \, ,
$$
while $M$ is a critical point of the functional {\rm (\ref{mc})} if, and
only if, its mean curvature function $h$ satisfies the equation
$$
\mbox{}\hfill 2\Delta h = 2cn h +2h\| \a\|^2 -h^3 \, . 
$$
\qed
\end{theorem}

\section{The Laplacian of the second fundamental form in space forms}

The Laplacian of the second fundamental form of the immersion
$f: M \rightarrow \tm$ 
was initially computed by Simons \cite{si} under the assumption of 
minimality on $f(M)$. Berard \cite{be} wrote down the general result, though
he applied it to immersions in space forms under the assumption that the mean 
curvature vector of the immersion $H$ was covariantly constant. In this 
section, we recall this fact in general. Once we show the minimality
of the critical points we consider of the functionals (\ref{sf}) and
(\ref{mc}), the result can be used in the same manner as before to
explain the interpretation of the gap theorem that we now make. 

Given any tensor field $Z$ on $M$, we set 
$\n^2_{X,Y}Z=(\n^g_X \n^g_Y -\n^g_{\n^g_X Y})Z$, an operator of order $2$. It
is tensorial in $X,Y$, and its symmetric properties are captured in the 
relation $(\n^2_{X,Y}-\n^2_{Y,X})Z=R^g(X,Y)Z$. The Laplacian is given by 
$\n^2=\sum_{i=1}^n (\n^g_{e_i}\n^g_{e_i}-\n^g_{\n^g_{e_i}e_i})$, where
$\{ e_i\}_{i=1}^n$ is a $g$-orthonormal frame on $M$. If we think
of the shape operator as the linear mapping $A: \nu(M) \rightarrow 
{\mathcal S}(M)$, where ${\mathcal S}(M)$ is the bundle of symmetric bilinear
maps on $TM$, then we obtain an adjoint $\mbox{}^{\dagger}\! A: {\mathcal S}(M)
\rightarrow \nu(M)$ defined by $\< \mbox{}^{\dagger}\! A s,N\>=\< A_N,s\>$.
Let $\widetilde{A}=\mbox{}^{\dagger }\! A \circ A$ and 
$\lat = \sum {\rm ad}A_{\nu_j}{\rm ad}A_{\nu_j}$,
where
$\{ \nu_j\}_{j=1}^q$ is an orthonormal basis of $\nu(M)$ at the point.
Finally, we have the curvature type operators $\overline{R}{'}_W$ and
$\overline{R}(A)^W$ given by
$$
\begin{array}{rcl}
\< {R^{\tg}}'_W(X), Y\> & = & \sum_i 
\< (\n^{\tg}_X R^{\tg})(e_i,Y)e_i,W\>+
\< (\n^{\tg}_{e_i}R^{\tg})(e_i,X)Y, W\> \, , \vspace{1mm} \\
\< R^{\tg}(A)^W(X),Y\> & = & \sum_i 2( \< R^{\tg}(e_i,Y)\a (X,e_i), W\>
 +\< R^{\tg}(e_i,X)\a (Y,e_i),W\>) \\
& & -\< A_N X, R^{\tg}(e_i,Y)e_i\>-\<A_W(Y),R^{\tg}(e_i,X)e_i\> \\
& & \< R^{\tg}(e_i,\a(X,Y))e_i,W\> -2\< A_W e_i,R^{\tg}(e_i,X)Y\> \, ,
\end{array}
$$
where $W$ is a normal vector at the point.  

\begin{theorem}\cite[Theorem 2]{be}
We have that
\begin{equation}\label{lap}
\begin{array}{rcl}
\< \n^2 \a (X,Y),W\> & = & -\< A\circ \widetilde{A}(W)(X),
Y\>-\<\lat \circ A_W (X),Y\> + \< R^{\tg}(A)^W(X),Y\>+ \\ & &  \< 
{R^{\tg}}'_W(X), Y\> + \< \n^2_{X,Y}H,W\>+\< R^{\tg}(H,X)Y,W\> \\
& & + \< A_W Y, A_H X\> \, .
\hfill \qed
\end{array}
\end{equation}
\end{theorem}

We now assume that $(\tm,\tg)$ has constant sectional curvature $c$.
Then we have that
\begin{equation}\label{id1}
R^{\tg}(X,Y)Z=c(\< Y, Z\> X -\<X,Z\>Y) \, ,
\end{equation}
that ${R^{\tg}}'_W=0$, and that $R^{\tg}(A)^W$ is given by
\begin{equation}\label{id2}
R^{\tg}(A)^W=cn\left( A_W - \frac{2}{n}\< H,W\> \BOne_{TM}\right) \,  . 
\end{equation}

The following result parallels \cite[Corollary 3]{be} but leaves in the term
involving the second derivatives on $H$ that appear in (\ref{lap}), as
we do not assume that the mean curvature vector is covariantly constant.

\begin{corollary}\label{c34}
If $(\tm,\tg)$ has constant sectional curvature $c$, then
$$
\n^2 \a = -(A\circ \widetilde{A}+ \lat \circ A)
+cn\left( A - \frac{1}{n}\< H, \, \cdot \, \> \BOne_{TM}\right)
+\n^2_{\, \cdot \, ,\, \cdot \,}H+ A\circ A_H \, .
$$
\end{corollary}

{\it Proof}.  We identify the terms in the right side of (\ref{lap}). The
first two yield $-(A\circ \widetilde{A}+ \lat \circ A)$. Now, by (\ref{id1}),
we have that $\< R^{\tg}(H,e_i)e_j,\a(e_i,e_j)\>= c\| H\|^2$, and so 
$R^{\tg}(H,\, \cdot \,)\, \cdot \, =c\< H, \, \cdot \, \> \BOne_{TM}$. This 
and (\ref{id2}) show that the third and and fifth terms in the right
side of (\ref{lap}) yield the next term in the right side of the 
stated expression for $\n^2 \a$. The fourth term in the right side
of (\ref{lap}) is zero, the next is $\n^2_{\, \cdot \, ,\, \cdot \,}H$, and
the last is $A\circ A_H$. 
\qed
 
We develop now a basic inequality for the analysis 
of (\ref{sf}) and (\ref{mc}), still in the case 
 when $(\tm,\tg)$ has constant sectional 
curvature $c$.  This inequality parallels and extends that of 
Simons \cite[Theorem 5.3.2]{si}, and reproduces Simons' result 
\cite[Theorem 5.3.2, Corollary 5.3.2]{si}
 when dealing with critical points of our functionals that are minimal 
submanifolds. As we will see, this embodies already the fact that the gap 
phenomena is one 
about the lowest critical values of the said functionals.

We recall Simons' \cite[Lemma 5.3.1]{si} inequality 
\begin{equation}\label{se}
\< A\circ \widetilde{A}+ \lat \circ A, A\> \leq \left( 2-\frac{1}{p}\right)
\| A\|^4 \, ,
\end{equation}
where $p$ is the codimension of $M$ in $\tm$. Then, by Corollary \ref{c34}, we
obtain the inequality
\begin{equation}\label{ine}
\begin{array}{rcl}
0 \leq {\displaystyle -\int \< \n^2 \a , \a\> d\mu } & = & 
{\displaystyle   \int \< A\circ \widetilde{A}+ \lat \circ A, A\> + 
\left( 2c\|H\|^2-cn \| \a \|^2\right) d\mu   
} \vspace{1mm} 
\\ & & {\displaystyle -\int  \left( \< \n^2_{e_i,e_j}H,\a(e_i,e_j)\>+
c\| H\|^2 + 
\< A_{\a(e_i,e_j)}e_j, A_H e_i\>\right)d\mu } \\ & \leq & 
{\displaystyle  \int \left( \left( 2-\frac{1}{p}\right)\| \a \|^4 +  
 c\| H\|^2 -cn\| \a \|^2 \right) d\mu  + } \vspace{1mm} 
\\ & & {\displaystyle -\int  \left( \< \n^2_{e_i,e_j}H,\a(e_i,e_j)\>+
\< A_{\a(e_i,e_j)}e_j, A_H e_i\>\right) }\, .
\end{array}
\end{equation}
Notice that when $H\equiv 0$, the formula above reduces to that of
Simons \cite[Theorem 5.3.2]{si}.

\section{Critical points of the Lagrangians with small density}

In this section we prove our main results. 

We consider, as above, the orthonormal frame
$\{ e_1, \ldots, e_n\}$ of tangent vectors to $(M,g)$ and the orthonormal
frame $\{ \nu_1, \ldots, \nu_p\}$ of the normal bundle $\nu(M)$ such that 
$H=h\nu_1$. We write $\a(e_i, e_j)=\sum_{r=1}^p h_{ij}^r \nu_r$. Then we have that
\begin{equation}\label{r1}
\text{$h=\sum_{i=1}^n h_{ii}^1$ and $\sum_{i=1}^nh_{ii}^r=0$ for
  $2\leq r\leq p$}\, , 
\end{equation} 

We denote by $\< \, \cdot \, , \, \cdot \, \>$ the standard metric on the space form
$S^{n+p}_c$. For immersions into it,   
the Ricci tensor and scalar curvature of $(M,g)$ are such that 
\begin{equation} \label{rtc}
r_g(e_j,e_k)=(n-1)c \< e_j,e_k\> +(hh_{jk}^1-
\sum_{i,r}h_{ij}^r h_{ik}^r ) \, ,
\end{equation} 
and
\begin{equation} \label{scc}
s_g = n(n-1)c + h^2- \| \a \|^2 \, ,
\end{equation}
respectively. 
Thus, under the hypothesis of Theorems \ref{mt} and \ref{mt2} we immediately derive that
$s_g \geq 0$ in the case when $c>0$.  

{\it Proof of Theorem \ref{mc}}. Suppose that $M\hookrightarrow \mb{S}^{n+p}$ is a 
critical point of $\Psi$ of constant mean curvature function. By Theorem \ref{th2}, 
we must have that 
\begin{equation}\label{in1}
0=2\int h\Delta h d\mu_g = \int h^2( 2n- 2\| \n_{e_i}^\nu \nu_1\|^2 -h^2 +
2 {\rm trace}\,A_{\nu_1}^2) d\mu_g  \, .
\end{equation}
Now the trace of $A_{\nu_1}^2$ is nonnegative and bounded above by $\| \alpha\|^2$. By
(\ref{es}) it follows that  
$$
\left(\frac{1}{2}-\lambda \right)
h^2 \leq n- \| \n_{e_i}^\nu \nu_1\|^2 -\frac{1}{2}h^2 +
 {\rm trace}\,A_{\nu_1}^2\leq  n- \| \n_{e_i}^\nu \nu_1\|^2 -\frac{1}{2}h^2 +
\| \alpha \|^2  
$$
is nonnegative and cannot be zero if $h\neq 0$. If $h=0$ then $\nabla^\nu \nu_1 =0$, and 
we have then that $0\leq \| \alpha \|^2 \leq n$. The desired result now follows by 
Theorem \ref{t1}. Indeed, the Clifford torus $\mb{S}^m(\sqrt{m/n})\times 
\mb{S}^{n-m}(\sqrt{(n-m)/n})$, $1\leq m < n$, have principal curvatures 
$\pm \sqrt{(n-m)/m}$ and $\mp \sqrt{m/(n-m)}$ with 
multiplicities $m$ and $n-m$, respectively, and by (\ref{rtc}), its Ricci curvature is 
bounded in between $n(m-1)/m$ and $n(n-m-1)/(n-m)$, while the real projective plane is 
embedded into $\mb{S}^4$ with an Einstein metric of scalar curvature $2/3$.   
\qed 
\medskip

We pause briefly to derive an elementary result to be used in the proof of Theorem
\ref{mt2}. This will unravel the geometric content of the constants that appear in 
the estimates (\ref{es2}), and the inequality of Theorem \ref{th6}.

\begin{lemma} \label{l11}
Let $(M,g)$ be a Riemannian manifold isometrically immersed into a background manifold 
$(\tm,\tg)$, and consider
the degree $1$-homogeneous function
$$
\frac{{\rm trace}\, \left( A_{\nu_1}\sum_j A_{\nu_j}^2\right)}{\| \alpha\|^2} \, .
$$
At a critical point we have that
$$
{\rm trace}\, \left( A_{\nu_1}\sum_j A_{\nu_j}^2\right)  = 
\frac{\| H\| ( \| \alpha \|^2 +2\| A_{\nu_1} \|^2 )}{n+2\|H\|^2/\| \alpha \|^2} \, ,
$$
the maximum occurs when $\| \alpha \|^2 = \| A_{\nu_1}\|^2$, and so
$$
{\rm trace}\, \left( A_{\nu_1}\sum_j A_{\nu_j}^2  \right) \leq \frac{3\| H\| \| \alpha \|^2}{2}\, .
$$
\end{lemma}

{\it Proof}. We have that
$$
{\rm trace}\, \left( A_{\nu_1}\sum_j A_{\nu_j}^2\right) = \sum_{k=1}^p \sum_{i,l,s=1}^n 
h_{is}^1 h_{sl}^k h_{li}^k \, ,
$$
and so the function of the $h_{ij}^k$s under consideration is defined by 
$$
\frac{{\rm trace}\, \left( A_{\nu_1}\sum_j A_{\nu_j}^2\right)}{\| \alpha\|^2}
 = \frac{\sum_{k=1}^p \sum_{i,l,s=1}^n h_{is}^1 h_{sl}^k h_{li}^k}{\sum_{k}\sum_{ij} (h_{ij}^k)^2}
$$
outside the origin, and extended by continuity everywhere. Its critical points subject to the 
constraints (\ref{r1}) are the solutions of the system of equations
$$
\sum_{k=1}^p \sum_{l=1}^n h_{ul}^k h_{lv}^k \delta_{1r} +\sum_{i=1}^n 
h_{ui}^1 h_{iv}^r + \sum_{s=1}^n h_{vs}^1 h_{su}^r =
({\sum_{k}\sum_{ij} (h_{ij}^k)^2})^2
\delta_{uv} \sum_{j=1}^p \lambda_j \delta_{jr} \, ,
$$
where $\lambda_1, \ldots, \lambda_p$ are the Lagrange multipliers. Here, 
$\delta $ is the Kronecker symbol.
If we multiply by $h_{uv}^r$ and add in $u,v$ and $r$, we obtain the relation
$$
\| \alpha \|^2
{\rm trace}\, \left( A_{\nu_1}\sum_j A_{\nu_j}^2\right) = \lambda_1 h \| \alpha \|^4 \, ,  
$$
while if we set $u=v$, $r=1$ and add in $u$, we obtain that
$$
\| \alpha \|^2(\| \alpha \|^2 +2 \| A_{\nu_1}^2 \|^2 )-2 h\,  
{\rm trace}\, \left( A_{\nu_1}\sum_j A_{\nu_j}^2\right) = \lambda_1 n \| \alpha \|^4 \, .  
$$
A simple algebraic manipulation yields the stated equality at critical points, and the statement
about the maximum is then clear. As we have assumed that $n\geq 2$, the inequality follows.
\qed
 
{\it Proof of Theorem \ref{mt2}}. Proceeding as above, by Theorem \ref{th2} 
we must have that 
\begin{equation}\label{in2}
0=2\int h\Delta h \, d\mu_g = \int h^2\left( 2-2 \| \n_{e_i}^\nu \nu_1\|^2 -\| \alpha \|^2 +
2\, {\rm trace}\,\frac{1}{h}A_{\nu_1} \sum_k A_{\nu_k}^2\right) d\mu_g  \, .
\end{equation} 
If estimates (\ref{es2}) hold for $\lambda \in [0,1)$, 
then $h=0$ and so $M$ is a minimal submanifold, and $\nabla^\nu \nu_1=0$. 
The desired result now follows by Theorem \ref{t1}. The argument is parallel to
the one already used in the proof of Theorem \ref{mt}. It suffices to add that the
symmetric Clifford torus $\mb{S}^m(\sqrt{1/2})\times \mb{S}^m(\sqrt{1/2})\subset
\mb{S}^{2m+1}$ is Einstein.   
\qed

{\it Proof of Corollary \ref{co4}}. There only needs to be observed why
these Riemannian manifolds are critical points of the total scalar curvature
under deformations of the metric $g$ in $M$ that can be realized by 
an isometric immersion of $M$ into $\mb{S}^{n+p}$. By (\ref{sc}), we have that
$$
\int s_g d\mu_g = n(n-1)\mu_g (M)+ \Psi(M)-\Pi(M) =\Theta(M)-\Pi(M)\, ,
$$
where $\mu_g(M)$ is the volume of the isometrically immersed manifold.
The critical manifold $M$ is separately a critical point of the functionals
$\mu_g(M)$, $\Psi(M)$ and $\Pi(M)$ under deformations of the immersion, 
hence a critical point of their linear combination above. 
\qed
\medskip

\begin{remark}
It is worth observing that the Klein bottle, a nonoriented manifold of zero 
Euler characteristic, can be embedded into $\mb{R}^4$ and therefore into
$\mb{S}^4$. But no embedding of such can be a critical point of $\Pi$ 
of constant mean curvature function and satisfying estimates (\ref{es2}).
\qed 
\end{remark}
\medskip

{\it Proof of Theorem \ref{th5}}. By Theorem \ref{th2}, we have that
\begin{equation}\label{inl}
0\leq 2\int h\Delta h \, d\mu_g = \int h^2(-2n- 2\| \n_{e_i}^\nu \nu_1\|^2 -h^2 +
2 {\rm trace}\,A_{\nu_1}^2) d\mu_g  \, .
\end{equation}
But $0\leq {\rm trace} \, A_{\nu_1}^2 \leq \| \alpha \|^2$, and so 
$$
-2n-h^2+2{\rm trace} \, A_{\nu_1}^2 \leq -2n-h^2+2\| \alpha \|^2\leq 0\, .
$$ 
It follows that either $h=0$, or that 
$$
-2n-h^2+2{\rm trace} \, A_{\nu_1}^2 = -2n-h^2+2\| \alpha \|^2 = 0\, ,
\quad \nabla_{e_i}^{\nu} \nu_1 =0\, .
$$
In the latter case, we have that $h_{ij}^r =0$ for all $r\geq 2$, and the vector $H=h\nu_1$ is 
a covariantly constant section of the normal bundle $\nu(M)$. The desired result follows.   
\qed
\medskip

{\it Proof of Theorem \ref{th6}}. We use once again the critical point equation given by
Theorem \ref{th2}, and obtain that  
$$
0\leq 2\int h\Delta h \, d\mu_g = \int h^2(-2-2\| \n_{e_i}^\nu \nu_1\|^2 -\| \alpha \|^2 +
2\frac{1}{h} {\rm trace}\,A_{\nu_1}\sum_j A_{\nu_j}^2) d\mu_g  \, .
$$ 
By Lemma \ref{l11}, we have that
$$
-2-\| \alpha \|^2 + 2\frac{1}{h} {\rm trace}\,A_{\nu_1}\sum_j A_{\nu_j}^2 
\leq -2 -\| \alpha \|^2+ 2\frac{ 3 \| \alpha \|^2 }{n+2 \| H \|^2 /\| 
\alpha \|^2}\, ,
$$
and the stated inequality is equivalent to the right side of this expression being nonpositive.
Thus, either $h=0$ or the right hand side of the expression above vanishes and the equality 
holds, and $\nabla^\nu \nu_1 =0$. In the latter case, $h_{ij}^r =0$ for all $r\geq 2$ and
$H=h\nu_1$ is a covariantly constant section of $\nu(M)$. 
\qed
\medskip

\noindent {\bf Acknowledgements.} Walcy Santos was partially supported by CNPq and 
Faperj of Brazil.


\begin{thebibliography}{69}
\bibitem{be}
P. B\'erard, {\it Simon's equation revisited}, Anais Acad. Brasil Ci\^encias
66 (1994) (4), pp. 397-403. 
\bibitem{bla}
W. Blaschke, {Vorlesungen \"uber Differentialgeometrie}, III. Springer
(1929), Berlin.
\bibitem{cdck}
S.S. Chern, M. do Carmo \& S. Kobayashi, {\it Minimal submanifolds of a 
sphere with second fundamental form of constant length}.  1970  
Functional Analysis and Related Fields (Proc. Conf. for M. Stone, Univ. 
Chicago, Chicago, Ill., 1968)  pp. 59-75 Springer, New York.
\bibitem{md}
M. Dajczer, {\it Submanifolds and Isometric Immersions}, Math. Lect. Ser. 13,
Publish or Perish, 1990.
\bibitem{glw}
Z. Guo, H. Li \& C. Wang, {\it The second variational formula for Willmore
submanifolds in ${\mathbb S}^n$}, Results in Math. 40 (2001), pp. 205-225.
\bibitem{bl}
H.B. Lawson Jr. {\it Local rigidity theorems for minimal hypersurfaces}, 
Ann. of Math. (2) 89 (1969), pp. 187-197.
\bibitem{hl}
H. Li, {\it Willmore hypersurfaces in a sphere}, Asian J. Math. 5 (2001), 
2, pp. 365-378.
\bibitem{s7}
S.R. Simanca, {\it The $L^2$-norm of the second fundamental form of isometric
immersions in Riemannian manifolds}, preprint 2013.
\bibitem{si}
J. Simons, {\it Minimal varieties in Riemannian manifolds}, Ann. Math. 2 
(1968), pp. 62-105. 
\bibitem{th}
G. Thomsen, {\it \"Uber Konforme Geometrie}, I. Grundlagen der Konformen 
Fl\"achentheorie, Abh. Math. Sem. Hamburg 3 (1923), pp. 31-56.
\bibitem{wi}
T.J. Willmore, Note on embedded surfaces.  An. \c Sti. Univ. ``Al. I. Cuza'' 
Ia\c si Sec\c t. I a Mat. (N.S.)  11B (1965), pp. 493-496. 
\end{thebibliography}
\end{document}